\documentclass[runningheads]{llncs}
\usepackage{graphicx}
\usepackage{amssymb}
\usepackage{url}
\usepackage{algorithm}
\usepackage{algpseudocode}
\usepackage{amsmath}
\DeclareMathOperator*{\avg}{avg}
\usepackage{color}

\begin{document}
\title{Node-Binded Communities for Interpolation on Graphs}


%
%
\author{Roberto Cavoretto\orcidID{0000-0001-6076-4115} \and
 Alessandra De Rossi\orcidID{0000-0003-1285-3820} \and
 Sandro Lancellotti\orcidID{0000-0003-4253-3561} \and
 Federico Romaniello\orcidID{0000-0003-1166-3179}}
\authorrunning{R. Cavoretto et al.}
%
\institute{Department of Mathematics \lq\lq Giuseppe Peano\rq\rq, University of Torino, via Carlo Alberto 10, 10123 Torino, Italy \\
\email{roberto.cavoretto@unito.it, alessandra.derossi@unito.it, sandro.lancellotti@unito.it, federico.romaniello@unito.it}}
\maketitle              
\begin{abstract}
Partition of unity methods (PUMs) on graphs represent strai\-ghtforward and remarkably adaptable auxiliary techniques for graph signal processing. By relying solely on the intrinsic graph structure, we propose the generation of a partition of unity through centrality measures and modularity. Subsequently, we integrate PUMs with a local graph basis function (GBF) approximation approach to achieve low-cost global interpolation schemes.

\end{abstract}
\section{Introduction}
Graph signal processing is a very popular tool when it comes to study the field of graph signals and, for this reason, several mathematical techniques such as filtering, compression, noise removal, sampling or decomposition are deeply investigated in the literature, see \cite{ort18,sta19}. Moreover, graphs naturally appear in a multitude of problems, as in social networks or traffic maps. However, in general, they have a complex structure that requires efficient and fast processing tools for their analysis.\\
Partition of unity methods (PUMs) allow to perform operations such as signal reconstruction from samples, classification of nodes, or signal filtering locally on smaller portions of the graph, and, then, rebuilding the global signal from the local ones.  In meshfree approximation, the combination of radial basis functions (RBFs) with PUMs yields significantly sparser system matrices in collocation or interpolation problems, and, therefore, a considerable speed-up of calculations \cite{cav20}. In what follows, we will focus on approximation methods based on generalized translates of a graph basis function (GBF), see \cite{erb1} and \cite{erb2} for more details.\\
In this article, we provide a more adaptive technique for the selection of the partitions on the graph than the one presented in \cite{cde}.  In particular, by using a process that automatically find an optimal number of subdomain based on their modularity and on the underlying graph structure, we do not need the desired number of subdomains as an input for the GBF-PUM algorithm.\\
In Section \ref{preli} we give a brief overview of graph theory and positive definite GBFs for signal approximation on graphs. A description of our algorithm and how it can be integrated to the GBF-PUM proposed in \cite{cde} is provided in Section \ref{tbsd}. In the final Section \ref{ne} we conclude with some numerical examples.

\section{Preliminaries}\label{preli}
\subsection{Graph theory}
In this article, we will use connected, simple and undirected graphs. Indeed, a graph $G$ is defined as $G=(V(G),E(G))$ where:
\begin{itemize}
    \item $V(G)$ is the \emph{vertex set} of the graph;
    \item $E(G) \subseteq V(G) \times V(G)$ is the \emph{edge set} of the graph containing all edges of the form $(u,v): u \neq v$. 
    \item $\textbf{A}$ is the \emph{adjacency matrix} of $G$, where $A_{ij}=1$ if $(v_i,v_j)$ is an edge, zero otherwise;
    \item A general symmetric graph Laplacian of the form $\textbf{L}= \textbf{D} - \textbf{A}$, where $\textbf{D}$ is the diagonal matrix with the degree of the vertices as entries.
\end{itemize}

A \emph{graph signal} $x: V(G) \rightarrow \mathbb{R}$ is a function on the vertices $V$ of $G$. We denote the $n$-dimensional vector space of all graph signals by $\mathcal{L} (G)$.
A graph $H=(V(H),E(H))$ is called a \emph{subgraph} of $G=(V(G),E(G))$ if and only if $V(H) \subseteq V(G)$ and $E(H) \subseteq E(G)$. A collection of disjoint subgraphs of $G$ which covers $V(G)$ is a \emph{partition} of $G$ and the subgraphs are called \emph{communities}. Given $u \in V(G)$, the set $N(u)=\left\{v\in V(G)| (u,v) \in E(G)\right\}$ is the \emph{neighbourhood} of $u$. Let $u,v \in V(G)$, the \emph{Jaccard similarity index} $J(u,v)= \frac{|N(u)\cap N(v)|}{|N(u)\cup N(v)|}$ measures the similarity between two vertices. It can easily be extended to communities $U,V$ of $G$ as follows: $$J(U,V)=\avg\limits_{u\in U, v \in V}\left\{J(u,v)\right\}.$$\\
The \emph{modularity} $Q$ of a graph measures the strength of division of a graph into communities $c_1,\dots, c_n$. Graphs with high modularity have dense connections between the vertices within communities but sparse connections between vertices in different communities. It is defined as 
$$Q=\frac{1}{2m}\sum_{i,j}\left(A_{ij}-\frac{k_ik_j}{2m}\right)\delta(c_i,c_j),$$
where $m$ is the number of edges,  $k_i$ is the degree of $i$ and $\delta(c_i,c_j)$ is 1 if $i$ and $j$ are in the same community, zero otherwise. \emph{Katz centrality} computes the relative influence of a vertex within a graph by measuring the number of the immediate neighbours and also all other vertices in the network that connect to the vertex under consideration through these immediate neighbours. Connections made with distant neighbours are, however, penalized by an attenuation factor $\alpha$, usually set as $0.5$. The Katz centrality for vertex $v_i$ is: 
$$C_{Katz}(v_i)= \sum_{k = 1}^{\infty}  \sum_{j = 1}^n  \alpha^k (A^k)_{ij}.$$

\subsection{Positive definite GBFs for signal approximation on graphs}
GBFs prove to be straightforward and effective instruments for the interpolation and approximation of graph signals, especially in scenarios where only a limited number of samples is available. Notably, employing positive definite GBFs tends to yield superior results. The theoretical underpinnings closely align with a parallel theory in scattered data approximation using positive definite RBFs in Euclidean space, as detailed in \cite{SW} and \cite{Wen}. Within this graph-based methodology, the approximation spaces are constructed by employing generalized shifts of a GBF, denoted as $f$, and are articulated in terms of a convolution operator. Numerical examples presented in Section \ref{ne} rely on polyharmonic splines on graphs, which are an example of a positive definite GBF based on the kernel:
$$\textbf{K}_{f_{(\epsilon\textbf{I}_n+\textbf{L})^{-s}}}=(\epsilon \textbf{I}_n+\textbf{L})^{-s}=\sum_{k=1}^{n} \frac{1}{(\epsilon + \lambda_k)^s}u_ku_k^T,$$
where $u_1,\dots,u_n$ are the eigenvectors of the Laplacian $\textbf{L}$ which form an orthonormal basis for the space $\mathcal{L}(G)$. 
If $\lambda_1$ is the largest eigenvalue of $\mathbf{L}$, this kernel is positive definite for $\epsilon > -\lambda_1$ and $s > 0$ (see, \cite{pes} and \cite{wanawa}).

\section{Topology Based Subdomain Detection}\label{tbsd}
Signal interpolation on graphs with PUM requires overlapping subdomains with the property that at least one vertex of the subdomain is an interpolation node. In what follows we describe how it is possible to find such subdomains for PUM on graphs, based on the structural topology of the underlying structure.
\subsection{Overlapping communities detection}
Communities detection is an edge-cutting topic in graph theory, indeed different approaches may be used, see \cite{newman} for a detailed overview. Our method is a divisive technique where the centrality of the interpolation nodes is the main criterion for splitting. In particular, we extend the procedure presented in \cite{gralg} allowing overlapping among communities. Our iterative algorithm repeats the following steps until some conditions on the presence of interpolation nodes and modularity are satisfied. It can be summarised as follows:
\begin{algorithm}
\caption{Node-Binded Communities Detection}
\hspace*{\algorithmicindent} \textbf{Input:} The graph $G$, the set of interpolation nodes $W$.
\begin{algorithmic}
\State Set $C=G$, $Q=-1$, $Q^{\prime}=-\frac{1}{2}$;
\While{$|C| \leq |W|$ and $Q^{\prime}>Q$:}
    \State For each $c$ in $C$  try to split and update $C$, $Q$ and $Q^{\prime}$;
\EndWhile
\State Merge small communities with the most similar big one and update $C$;
\State For each $c$ in $C$ expand $c$ to create overlapping subdomain.
\end{algorithmic}
\hspace*{\algorithmicindent} \textbf{Output:} The set of overlapping communities $C$. 
\end{algorithm}


In detail, in the splitting process we find the two interpolation nodes with the highest Katz centrality in each community, if they exist, and divide the community into two sub-components, each one containing one of the samples. When the splitting process terminates, the \emph{small} communities containing less than $2\%$ of vertices are merged with the most similar big community, according to the Jaccard index. Finally, the overlapping expansion is performed for each community, where the more edges a community shares with the others, the more it is augmented. 

\subsection{GBF-PUM approximation on graphs}
Partition of unity is a widespread tool used in mesh-free approximation problems as it sharply reduces the computational cost when related subproblems are solved \cite{fas15}. In \cite{cav22}  and \cite{cde} this technique has been adapted to graph signal interpolation and it may be applied to social networks or traffic map problems because of their underlying graph structure. The algorithm is made of the following steps:
\begin{algorithm}
\caption{GBF-PUM Approximation on Graphs}
\hspace*{\algorithmicindent} \textbf{Input:} The graph $G$, the set of communities $C$, the set of interpolation nodes $W$.
\begin{algorithmic}
\State Construct a partition of unity subordinate to the communities in $C$;
\State For all subgraphs $c$ in $C$ calculate the local Laplacian $\textbf{L}^c$;
\State For all subgraphs $c$ in $C$ construct the local GBF kernel;
\State For all subgraphs $c$ in $C$ calculate the local GBF approximant;
\State Create a global GBF-PUM approximation from the local ones.        
\end{algorithmic}
\hspace*{\algorithmicindent} \textbf{Output:} GBF-PUM approximant of the signal.
\end{algorithm}


We note that overlapping subdomains are needed for a good PUM. In addition to this, we remark that the global approximation is a weighted sum of the local approximants by means of weighted functions defined by each subdomain. Further details on the error of the global approximant can be found in \cite{cde}.

\section{Numerical Examples}\label{ne}
In this section our framework is applied on the Minnesota Road Graph \cite{minnroadgraph}, given a signal on its vertex set. This graph has 2642 vertices and 3304 edges; we test the algorithm on sets of 200, 400, 600 and 800 sampling nodes used as interpolation nodes. We remark here we do not require the number of communities as an input and this is the main difference with respect to \cite{cav22} and \cite{cde}.\\

\begin{figure}[h]
    \centering
    \includegraphics[width=.5\textwidth]{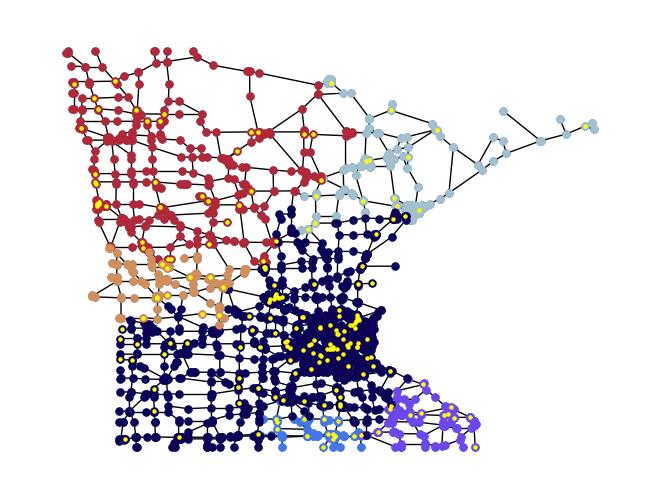} \\
\includegraphics[width=.5\textwidth]{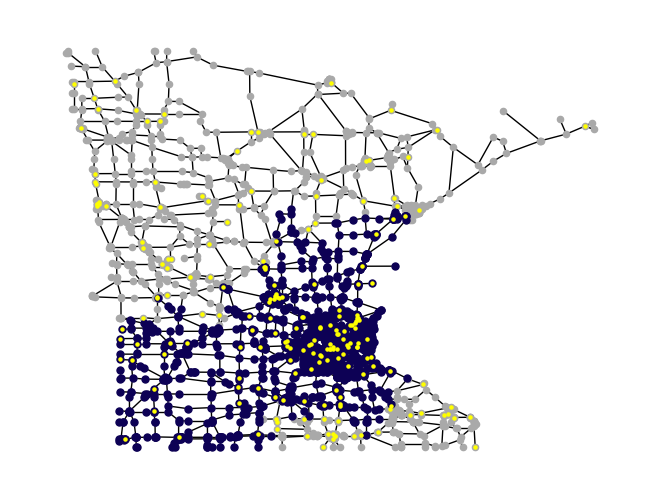}\hfill
\includegraphics[width=.5\textwidth]{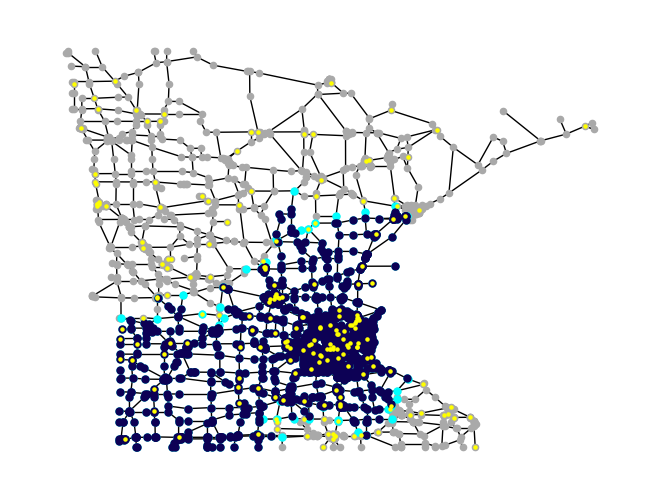}
    \caption{On top the partition of the vertices in 6 different disjoint subdomains is shown. It is obtained using 200 interpolation points, marked in yellow. On bottom a subdomain is highlighted to show how it is expanded to create the overlap (left to right).}
    \label{subdomain_extension}
\end{figure}

\begin{table}[h]
    \centering
    \begin{tabular}{|c|c|c|c|}
    \hline
    N       & Communities &     RRMSE &     time (s) \\
    \hline
    200     & 6           &    1.102e-01 & 1.18e02 \\
    \hline
    400     &  11         &  7.047e-02 & 1.26e02 \\
    \hline
    600     &  6          &  5.456e-02 & 1.26e02 \\
    \hline
    800     &  9            &  4.499e-02 & 1.69e02  \\
    \hline
    \end{tabular}
    \caption{A summarisation of the numerical results for increasing number of sampling nodes. The number of communities, RRMSE and computational time are also shown.}
    \label{tabellaprove}
\end{table}

It is worth noting that the numerical results on the Relative Root Mean Square Error (RRMSE) are of the same order as the ones obtained by the authors in \cite{cde} with 6 subdomains and augmentation of communities with neighbours at distance 2, given as input. We observe here that the number of communities in our algorithm depends only on the centrality of the interpolation nodes and it is not given as input. Moreover, our criterion for the augmentation of the communities depends on the ratio $r(v)$ between the number of neighbours of a vertex $v$ inside its community and the total number of its neighbours:
\begin{enumerate}
    \item if $r(v)\leq0.4$: expand with neighbours at distance 2;
    \item if $0.4<r(v) \leq 0.8$: expand with neighbours at distance 1.
\end{enumerate}
The numerical results are summarised in Table \ref{tabellaprove}. Figure \ref{subdomain_extension} highlights the 6 communities found by our algorithm with the 200 sampling nodes marked in yellow (top), then we show how one of them is augmented to create the overlapping (bottom).\\

 All tests were carried out on the infrastructure for high-performance computing \emph{MathHPC}, virtual cloud server of the structure \emph{HPC4AI} (High-Performance Computing for Artificial Intelligence: \url{https://www.dipmatematica.unito.it/do/progetti.pl/View?doc=Laboratori_di_ricerca.html}).

\subsection*{Acknowledgments}  

The authors sincerely thank the reviewers for the careful reading and valuable comments on the paper. This research has been accomplished within the RITA \lq\lq Research ITalian network on Approximation\rq\rq\ and the UMI Group TAA\lq\lq Approximation Theory and Applications\rq\rq. This work has been supported by the INdAM--GNCS 2022 Project \lq\lq Computational methods for kernel-based approximation and its applications\rq\rq, code CUP$\_$E55F22000270001, and by the Spoke ``FutureHPC \& BigData'' of the ICSC – National Research Center in "High-Performance Computing, Big Data and Quantum Computing", funded by European Union – NextGenerationEU. Moreover, the work has been supported by the Fondazione CRT, project 2022 \lq\lq Modelli matematici e algoritmi predittivi di intelligenza artificiale per la mobilit$\grave{\text{a}}$ sostenibile\rq\rq.


%
%
%
%

\end{document}